\documentclass[11pt]{article}
\usepackage{cite}
\usepackage{amsmath,amssymb,amsfonts}
\usepackage[margin=1in]{geometry}
\usepackage{textcomp}
\usepackage{amssymb}
\usepackage{graphicx}
\usepackage{bm}
\usepackage{algorithm2e}
\usepackage[noend]{algpseudocode}
\usepackage{subfig}
\usepackage[dvipsnames]{xcolor}
\usepackage{hyperref}
\newenvironment{spmatrix}[1]
{\def\mysubscript{#1}\mathop\bgroup\begin{pmatrix}}
	{\end{pmatrix}\egroup_{\textstyle\mathstrut\mysubscript}}

\newcommand{\eq}[2]{\begin{equation}\label{#1} #2 \end{equation}}
\newcommand{\al}[1]{\begin{align*} #1 \end{align*}}

\newcommand{\rf}[1]{(\ref{#1})}
\newcommand{\eqs}[2]{\begin{equation} \label{#1}\begin{split} #2\end{split} \end{equation}}

\newcommand{\norm}[1]{\left\lVert#1\right\rVert}

\renewcommand{\a}{{\alpha}}

\newcommand{\ve}{\varepsilon}
\newcommand{\tr}{\triangledown}
\newcommand{\td}{\tilde}

\newcommand{\g}{\gamma }

\newcommand{\s}{\sigma }

\renewcommand{\k}{\kappa}

\renewcommand{\d}{\delta}
\renewcommand{\t}{\theta}

\renewcommand{\l}{{\lambda}}

\newcommand{\E}{\mathbb{E}}

\renewcommand{\it}{_{i,t}}

\newcommand{\trs}{^\intercal}
\newcommand{\ik}{_{i,k}}
\newcommand{\jk}{_{j,k}}

\def\*#1{\bm{#1}}
\def\@#1{\mathcal{#1}}

\def\BibTeX{{\rm B\kern-.05em{\sc i\kern-.025em b}\kern-.08em
		T\kern-.1667em\lower.7ex\hbox{E}\kern-.125emX}}
\providecommand{\keywords}[1]{\textbf{\textit{Keywords:}} #1}

\begin{document}
	\title{Distributed Networked Real-time Learning}
	\author{Alfredo Garcia\thanks{This work was supported by NSF ECCS-1933878 Award  and  Grant AFOSR-15RT0767.  },
		Luochao Wang, Jeff Huang, and Lingzhou Hong
	\thanks{ \{alfredo.garcia, wangluochao, jeffhuang,hlz\}@tamu.edu}}
	\date{Texas A\&M University}
	
	\maketitle
	
	\begin{abstract}
		Many machine learning algorithms have been developed under the
		assumption that data sets are already available in batch form. Yet in many
		application domains data is only available sequentially overtime via compute nodes in different geographic locations. 
		In this paper, we consider the problem of learning a model when streaming data cannot be transferred to a single location in a timely fashion. In such cases, a distributed architecture for learning relying on a network of interconnected ``local" nodes is required. 
		We propose a distributed scheme in which every local node implements stochastic gradient updates based upon a local data stream. To ensure robust estimation, a network regularization penalty is
		used to maintain a measure of cohesion in the ensemble of models.
		We show the ensemble average approximates a stationary point and  characterize the degree to which individual models differ from the ensemble average. We compare the results with federated learning to conclude the proposed approach is more robust to heterogeneity in data streams (data rates and estimation quality).
		We illustrate the
		results with an application to image classification with a deep learning model based upon convolutional neural
		networks.
	\end{abstract}

 \begin{minipage}{12cm}%
 	\centering
 		\textbf{keywords:} asynchronous computing, distributed computing, networks,
 	                   non-convex optimization, real-time machine learning. 
 \end{minipage}

\section{Introduction}

\label{intro} Streaming data sets are pervasive in certain application domains often involving a network of
compute nodes located in different geographic locations. However, most machine
learning algorithms have been developed under the assumption that data sets
are already available in batch form. When the data is obtained through a
network of heterogeneous compute nodes, assembling a diverse batch of data
points in a central processing location to update a model may imply significant
latency.  Recently, an architecture referred to as {\em federated} learning (FL, see e.g. \cite{FL6, FL3}) with a central server in proximity to local nodes has been proposed. In FL, each node implements updates to a machine learning model that is kept in the central server. This allows collaborative learning while keeping all the training data on nodes rather than in the cloud. In general, schemes that avoid the need to rely on the cloud for data storage and/or computation are referred to as ``edge computing".

With high data payloads, such architecture for real-time learning is subject to an 
\emph{accuracy} vs \emph{speed} trade-off due to asymmetries in data quality
vs. data rates as we explain in what follows.   

Consider nodes $i\in\{1,\dots,N\}$ generating data points $%
(x_{i,n},y_{i,n}),~n \in \mathbb{N}^+$ at different rates $\mu_i>0$ which are used for the instantaneous
computation of model updates $\theta_{k}, ~ k \in \mathbb{N}^+$ (striving to
minimize loss $\ell$). This setting could correspond for example with supervised deep
learning in real-time wherein gradient estimates (with noise variance $%
\sigma^2_i>0$) are computed via back-propagation in a relatively fast
fashion. Without complete information on $\sigma^2_i>0$, updating the model
parameters based upon every incoming data point yields high speed but possibly at the expense of low accuracy. For example, if the nodes producing noisier estimates
are also \emph{faster} at producing data, it is highly unlikely that an
accurate model will be identified at all. 

To illustrate this scenario, in Figure 1, we depict the performance of FL for {\em deep} convolutional neural networks with the MNIST dataset. In these simulations, each one of $N=5$ nodes sends data according to independent Poisson processes with $\mu_0=8$ and $\mu_i=1$, $i \in \{1, \dots, 4\}$. The fastest node computes gradient estimates based upon a {\em single} image whereas the slower nodes compute gradient estimates based upon a batch of $64$ images.

\begin{figure}[tb]
\centering
\vspace{-10pt}
\includegraphics[width=.5\textwidth,height=!]{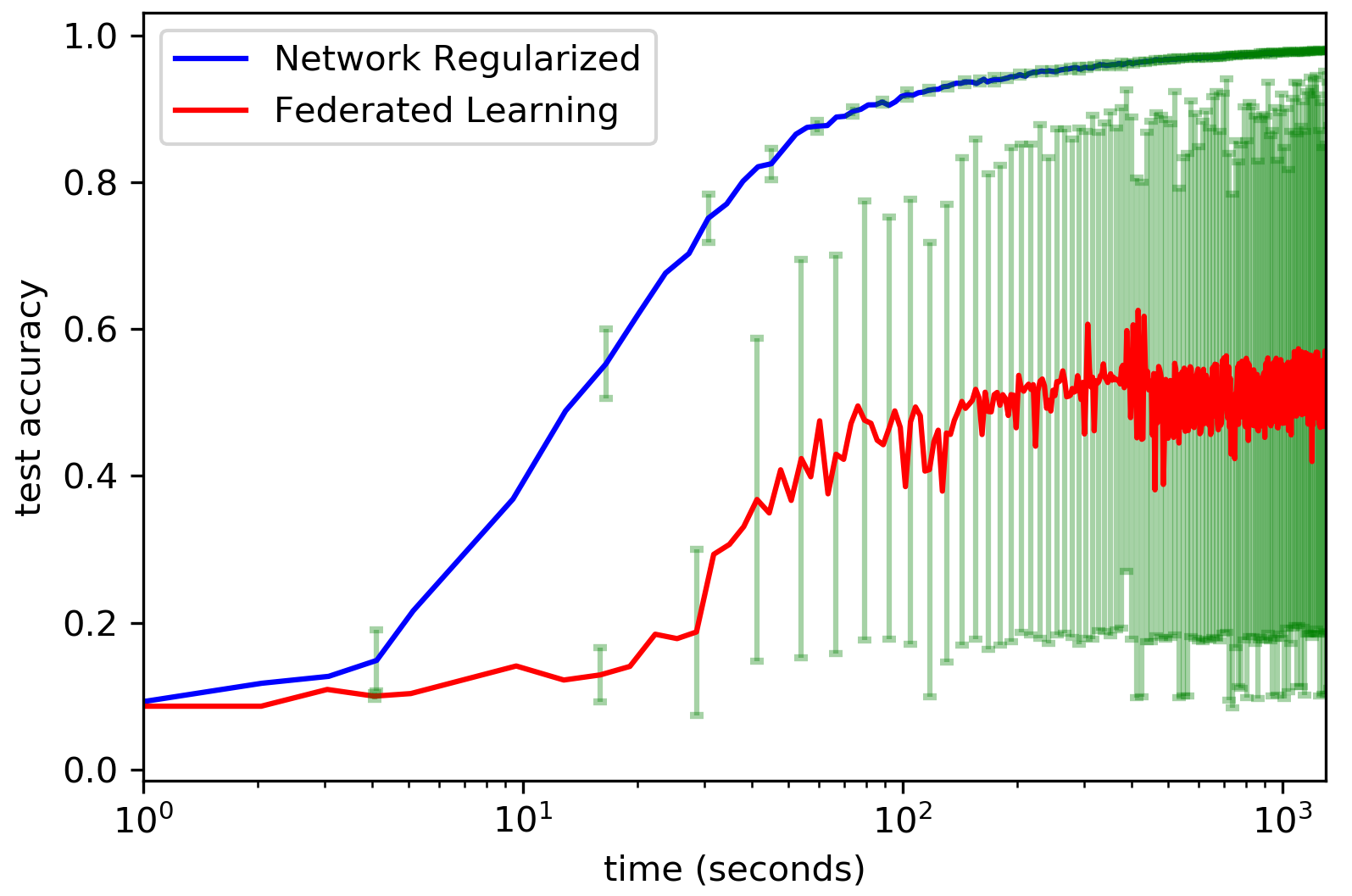}
\vspace{-7pt}
\caption{ Performance comparison for Deep Learning on MNIST with learning rate $\gamma=0.01$. The $95\%$ percentile is depicted with green lines.}
\label{fig1}
\end{figure}

This trade-off between \emph{speed} and \emph{%
precision} is mitigated in a distributed approach to real-time
learning subject to a network regularization penalty. In such an approach, each one
of $N>1 $ local nodes independently produce parameter updates based upon
a single (locally obtained) data point which speeds up computation.
Evidently, with increased noise, such a scheme may fail to enable the
identification of a reasonably accurate model. However, by adding a network
regularization penalty (which is
computed locally) a form of coordination between multiple local nodes is
induced so that the ensemble average solution is robust to noise.\footnote{Similar network regularization methods have been used in multi-task learning to account for inherent network structure in data sets (see e.g.\cite{Hallac2015, Sayed_1, Sayed_2, Sayed_3, Sayed_4,Sayed_5}). See section D.} Specifically, we show that the ensemble average solution approximates a stationary point and  that the approximation quality is
$\mathcal{O}(\frac{\sum_{i=1}^{N}\sigma^{2}_i}{N^{2}})$, which compares quite favorably with FL, which is highly sensitive to fast and inaccurate data streams.
We illustrate the results with an application to deep learning
with convolutional neural networks.

The structure of the paper is as follows. In section 2 we introduce the
distributed scheme  that combines stochastic gradient descent with network
regularization (NR). In section 3 we analyze the scheme and show that it
converges (in a certain sense) to a stationary point, we also compare its
performance with FL.
Finally, in section 4 we report the results from a testbed on deep learning
application to image processing, and in section 5 we offer conclusions.

\section{ A Network Regularized Approach to Real-time Learning}

\label{sec:network}

\subsection{Setup}

\label{sec:1}

We consider a setting in which data is made available sequentially overtime via nodes $i \in \{1,\dots, N\}$ in different geographic locations. We denote the $i$-th stream by $\{(\*x_{i,n},y_{i,n}):n\in \mathbb{N}^+\}$ and assume these data points are independent samples from a joint distribution $\mathcal{P}_i$. 

We also assume the data streams are independent but heterogeneous, i.e. $\mathcal{P}_i\neq\mathcal{P}_j, i \neq j$.
Each node strives to find
a parameter  specification $\theta \in \Theta \subset \mathbb{R}^p$ that
minimizes the performance criteria $\mathbb{E}%
_{\mathcal{P}_i}[\mathcal{L}(\*x_i,y_i;\theta)]$, where the loss function $\mathcal{L}(\cdot) \geq 0$ is continuously differentiable with respect to $\t$.  
Though data is distributed and heterogeneous, we consider a setting in which nodes {\em agree} on a common learning task.
This is formalized in the first standing assumption. 
Let $g_i(\theta)\triangleq \nabla_{\theta} \mathcal{L}(\*x_{i},y_{i};\theta)$ denote the gradient evaluated at $(\*x_{i},y_{i})\sim \mathcal{P}_i$, and assume $g_i(\t)$ is uniformly integrable.

\medskip 
\parindent 0cm 
\textbf{Assumption 0:} {\em For all $\theta \in \Theta$, and $i \in \{1,\dots,N\}$}:
\[
 \mathbb{E}%
_{\mathcal{P}_i}[g_i(\theta)]= \mathbb{E}%
_{\mathcal{P}_j}[g_j(\theta)].
\]
Let $\ell (\theta)$ denote the (ensemble) average expected loss:
\[\ell (\theta) \triangleq \frac{1}{N}\sum_{i=1}^{N} \mathbb{E}%
_{\mathcal{P}_i}[\mathcal{L}(\*x_i,y_i;\theta)].\]

By uniform integrability, $\nabla_{\theta} \E_{\@P_i}\@L(\*x_i,y_i;\t)= \E_{\@P_i}g_i(\t)$. Assumption 0 thus implies that $\E_{\@P_i}[g_i(\t)]=\tr \ell(\t)$ for all $i$ and $\t$.

 Let $\varepsilon _{i}(\theta) \triangleq g_i(\theta)-\tr \ell (\theta)$, then it holds that $\mathbb{E}[\ve_i(\theta)]=0$. We further assume:
 
\medskip \parindent 0cm \textbf{Assumption 1:} \emph{For all $\theta_i \in \Theta$, the random variables
$\{\varepsilon_{i}(\theta_i): i \in \{1,\dots, N\}\}$ are independent and} 
\[
 \mathbb{E}[\norm{\ve_i(\theta_i)}^2]\leq \sigma_i^{2}.
\]
Define $\s^2=\sum_{i=1}^N\s_i^2$. By independence of data streams:
\[
\mathbb{E}[\varepsilon _{i}(\theta)\trs\varepsilon
_{j}(\theta)]=\mathbb{E}[\varepsilon _{i}(\theta)]\trs\mathbb{E}[\varepsilon _{j}(\theta)]=0,
\]
for all $\t\in \Theta, j\in \left.{\{1,\dots,N\}}\right\backslash {\{i}\}$.

Streams generate data over time according to independent Poisson processes { $D_i(t)$ with rate $\mu _{i}>0$ and $D_i(0)=1$.}
We assume the time required to compute gradient estimates and/or exchange parameters locally among neighbors or with the central server are negligible compared to the time in between model updates. 
In what follows we make use of a virtual clock that produces ticks
according to an aggregate counting process $D(t)=\sum_{i=1}^{N}D_{i}(t)$ with
rate $\mu =\sum_{i=1}^{N}\mu _{i}$. Let $k\in \mathbb{N}^{+}$ denote the index set of ticks associated with the aggregate process. { Since we assume the parameter is updated once a data point arrives,  the $k$-th iteration is completed at the $k$-th tick. Index $k$ denotes the $k$-th step in the schemes described below. }

\subsection{Federated Real-time Learning}

In FL, gradient estimates are communicated to a central server where a model is updated as follows: 
\begin{equation}
 \theta _{k+1}=\t_k-\gamma \sum_{i=1}^{N}\mathbf{1}_{i,k}g_i(\theta_k),  \label{SG}
\end{equation}%
where $\g$ is the learning
rate, $\mathbf{1}\ik$ is a indicator of whether node $i$ performs the an update:  $\*1\ik=1$ if the next gradient estimate comes from the $i$-th
stream and $\mathbf{1}_{i,k}=0$ otherwise.

The algorithmic scheme described in (\ref{SG}) was first analyzed in 
\cite{Wright} for data in batch form and has been used in the recent literature on asynchronous
parallel optimization algorithms (see for example \cite{Liu}, \cite{Lian}
and \cite{Duchi_async}). As Figure 1 suggests,
with heterogeneous data streams, the scheme in (\ref{SG})
trades off speed in producing parameter updates at the expense of
heterogeneous noise in gradient estimates. In what
follows we introduce a distributed approach that relies on a network
regularization penalty to ensure the ensemble average approximates a
stationary point (i.e. a choice of parameters with null gradient). We will
show that in such a networked approach the trade-off between precision and
speed is mitigated.

\subsection{A Distributed Approach with Network Regularization}

In the NR scheme, we consider a network of local compute nodes which we model as a graph $\mathcal{G%
}=(\mathcal{N},\mathcal{E})$, where $\mathcal{N}=\{1,\dots ,N\}$ stands for
the set of nodes and $\mathcal{E}\subseteq \mathcal{N}%
\times \mathcal{N}$ is the set of links connecting nodes. Let $A=[\alpha
_{ij}]\in \mathbb{R}^{N\times N}$ be the adjacency matrix of $\mathcal{G}$,
where $\alpha _{ij}\in \{0,1\}$ indicates whether node $i$ communicates with node $j$: $\a_{ij}=1$ if two nodes can exchange  local information and $\a_{ij}=0$ otherwise.

In this scheme, each local node $i$ performs model updates according to a linear combination of local gradient estimate and the gradient of a consensus potential:
$$\mathcal{F}(\*{\t} )=\frac{1}{4}\sum_{i}\sum_{j\neq i}\alpha	_{ij}\Vert \theta _{i}-\theta _{j}\Vert ^{2},$$ 
where 
 $\*{\t} _{t}\trs=(\theta _{1,t}\trs,\dots, \theta _{N,t}\trs)\in \mathbb{R}^{p\times N}$.
 The consensus potential is a measure of similarity across local models.\footnote{This consensus potential has been used in the literature of opinion dynamics (see e.g. \cite{Friedkin}).} 
  The  update performed by node $i$ is of the form: 
\begin{equation}
\theta _{i,k+1}=\theta _{i,k}-\gamma\mathbf{1}_{i,k}[g_i(\theta
_{i,k})+ a\tr\@F\ik],
\label{flocking}
\end{equation}
where 
$a>0$ is a regularization parameter, and
\[
\tr \@F\ik \triangleq \nabla_{\theta_i}\@F(\*{\t}_k)=\sum_{j \neq i}\a_{ij}(\t\ik-\t\jk).
\]

Note that the basic iterate (\ref{flocking}) can be interpreted as a stochastic gradient approach to solve the local problem: 
\[
\min_{\theta_i}[\mathbb{E}_{\mathcal{P}_i}(\theta_i) +a
\mathcal{F}(\boldsymbol{\theta})],
\]
in which the objective function is a linear combination of loss and consensus potential.\footnote{{This interpretation is not novel (see e.g. \cite{Passino}, \cite{Pu} for its use in swarm (flocking) optimization and in multi-task learning \cite{Sayed_1, Sayed_2, Sayed_3, Sayed_4,Sayed_5}).}} 
When $a=0$, each local node ignores the neighboring models.  For large values of $a>0$, the resulting dynamics reflect the countervailing effects of seeking to minimize consensus potential and improving model fit. With highly dissimilar initial models, each local node largely ignores its own data and opts for updates that lead to a model that is similar to the local average. Once approximate consensus is achieved, local gradient estimates begin to dictate the dynamics of model updates.

In what follows it will be convenient to rewrite (\ref{flocking}) as follows: 
\eq{eq:basic flocking_2}{
\theta _{i,k+1}=\theta _{i,k}-\gamma\mathbf{1}_{i,k}\left[ \tr \ell
(\theta _{i,k})+a\tr \mathcal{F}_{i,k}+\varepsilon _{i,k}\right].}
Given that local nodes independently update and maintain
their own parameters, the network regularized scheme is not subject to the
possibility of biased gradient estimates stemming from update delays in FL (see
\cite{Leblond}).

\subsection{Literature Review}
The scheme proposed in (\ref{flocking}) has already been considered in the machine learning literature. In a series of papers (see \cite{Sayed_1, Sayed_2, Sayed_3, Sayed_4, Sayed_5})), the authors consider an approach to {\em multi-task} learning based upon a network regularization penalty as in (\ref{flocking}). This paper focuses on distributed {\em single-task} learning. In contrast to the papers referred above, we consider a non-convex setting with heterogeneous nodes asynchronously updating their respective models at different rates over time.

{ The scheme proposed in (\ref{flocking}) is also related to the literature on consensus optimization (see e.g. \cite{nedic2009distributed}, \cite{Lian}, \cite{Yin}). However, the proposed approach can not be interpreted as being based upon {\em averaging over local} models as in consensus-based optimization.} In that literature, the basic iteration is of the form:
\begin{equation*}
\theta _{i,k+1}=\sum_{j} W_{i,j,k}\theta _{j,k}-\gamma g
(\theta _{i,k}),
\end{equation*}
where ${\bf W}_k \in \mathbb{R}^{N \times N}$ is doubly stochastic and $g
(\theta _{i,k})$ is a noisy gradient estimate. Indeed one can rewrite (\ref{flocking}) as:
\begin{equation*}
\theta _{i,k+1}=\sum_jW_{i,j}\theta_{j,k}
- \gamma \mathbf{1}_{i,k}g_{i}(\theta_{i,k}),
\label{network}
\end{equation*}
with $W_{i,i}=1-\gamma a\sum_j \alpha_{i,j}$ and $W_{i,j}=\gamma a\sum_j \alpha_{i,j}$. However, the resulting matrix ${\bf W}$ is {\em not doubly stochastic} in general since we only require $a>0$. Thus, the approach to consensus in (2) can not be interpreted as being based upon {\em averaging over local} models as in consensus-optimization.

{The algorithms proposed in \cite{Lian} and \cite{Yin} are designed for {\em batch} data while our approach deals with {\em streaming} data. For example, in \cite{Lian}, each node uses the {\em same} mini-batch size for estimating gradients while in our approach gradient estimation noise is {\em heterogeneous}}. In addition, the algorithms proposed in \cite{Lian} and \cite{Yin}, every node is {\em equally likely} to be selected at each iteration to update its local model. 
In contrast, in our approach data streams are heterogeneous so that certain nodes are {\em more likely} to update their models at any given time. {Finally, in \cite{Lian} the objective function (loss)
is defined with respect to a distribution that is biased towards the
nodes that update more often. This is in contrast to the objective
function defined in this paper (i.e. $\ell(\theta)$), where every node
contributes to the global distribution with the same weight regardless of their updating frequency.}

\section{Analysis}

\label{statements} In this section, we show the NR scheme
converges (in a certain sense) to a stationary point. To that end we study stochastic processes $\{\theta _{i,k}:k>0\}$ associated with each one of
the $N>1$ nodes in the network regularized approach. The proofs are given
in the appendix. 
We make the following standing assumptions:

\medskip \parindent 0cm \textbf{Assumption 2}: \emph{The graph $\mathcal{G}$
corresponding to the network of nodes is undirected ($A=A\trs$) and
connected, i.e., there is a path between every pair of vertices.}

\medskip \parindent 0cm \textbf{Assumption 3}
\label{asp:Lipschitz} \emph{%
(Lipschitz) $\norm{\tr \ell(\theta )-\tr \ell(\theta')}\leq L\norm{\theta -\theta'} $  for some $L>0$ and for all $%
\theta ,\theta ^{\prime }$}. \label{sec: main}

\subsection{Preliminaries}

The ensemble average $\bar{\theta}_{k}\triangleq\frac{1}{N}%
\sum\nolimits_{i=1}^{N}\theta _{i,k}$ plays an important role in
characterizing the performance of the network regularized scheme. To this
end, we analyze the process $\{\overline{V}_{k}:k>0\}$ defined as 
\begin{equation*}
\overline{V}_{k}\triangleq \frac{1}{N}\sum_{i=1}^{N}\Vert \theta _{i,k}-\bar{%
\theta}_{k}\Vert ^{2}.
\end{equation*}
Let $e_{i,k}\triangleq \theta _{i,k}-\bar{\theta}_{k}$ and $%
V_{i,k}\triangleq\Vert e_{i,k}\Vert ^{2}$, then $\overline{V}_{k}=\frac{1}{N}%
\sum_{i=1}^{N}V_{i,k}$. We now introduce some additional notations. Let $\deg (i)$ denote the degree
	of vertex $i$ in graph $\mathcal{G}$ and  $\overline{d}:=\max_{i}\deg (i)$%
	. Let $\mathbb{E}[\overline{V}_{k+1}|\*{\t}_k]$ denote the
	conditional expectation of $\overline{V}_{k+1}$ given $\mathbf{\theta }_{k}$%
	. We define $\mu _{\max } =\max\{\mu_i:1\leq i \leq N\}$ and $\mu _{\min }
	=\min\{\mu_i:1\leq i \leq N\}$.   
We first prove two intermediate results.

\medskip \textbf{Lemma 1} \emph{Suppose Assumptions 0, 1 and 2 hold. It holds
that} \vspace{-0.25cm} 
\begin{align*}
\overline{V}_{k+1}& =\overline{V}_{k}-\frac{2}{N}\sum_{i=1}^{N}\gamma e_{i,k}\trs\mathbf{1}_{i,k}\left[ \tr \ell (\theta _{i,k})+a\tr 
\mathcal{F}_{i,k}\right]  \\
& \qquad -\frac{2}{N}\sum_{i=1}^{N}\gamma e_{i,k}\trs\varepsilon _{i,k}%
\mathbf{1}_{i,k}+\frac{1}{N}\sum_{i=1}^{N}\gamma^{2}\Vert \delta
_{i,k}\Vert ^{2},
\end{align*}%
\emph{where} $\delta _{i,k}=\delta _{i,k}^{f}+\delta _{i,k}^{g}+\delta
_{i,k}^{n}$, \emph{and} 
\begin{align*}
\delta _{i,k}^{f}\triangleq \tr \ell (\theta _{i,k})\mathbf{1}_{i,k}-
\tr\bar{\ell} _{k},& ~~~\tr\bar{\ell}_{k}\triangleq \frac{1}{N}%
\sum_{j=1}^{N}\tr \ell (\theta _{j,k})\mathbf{1}_{j,k}, \\
\delta _{i,k}^{g}\triangleq a(\tr \mathcal{F}_{i,k}\mathbf{1}_{i,k}-\tr\bar{\mathcal{F}}_{k}),& ~~~\delta _{i,k}^{n}\triangleq \varepsilon _{i,k}%
\mathbf{1}_{i,k}-\frac{1}{N}\sum_{j=1}^{N}\varepsilon _{j,k}\mathbf{1}_{j,k},
\\
~\tr\bar{\mathcal{F}}_{k}& \triangleq \frac{1}{N}\sum_{j=1}^{N}\tr 
\mathcal{F}_{j,k}\mathbf{1}_{j,k}. \\
&
\end{align*}%
\vspace{-0.25cm} \label{delta^fgn}

\medskip \textbf{Lemma 2} 
\label{lem2} 
{\em Suppose Assumptions 0, 1, 2 and 3 hold.  
Let $\xi={\mu _{\max }}/{\mu _{\min }}$, then:
\begin{equation*}
\mathbb{E}[\overline{V}_{k+1}|\*{\t}_k]\leq (1+\frac{\kappa \gamma }{N})%
\overline{V}_{k}+\frac{4\gamma^{2}\xi}{N}\left\Vert \tr \ell (\bar{%
\theta}_{k})\right\Vert ^{2}+\frac{\gamma^{2}\xi\sigma ^{2}}{N^{2}},
\end{equation*}%
where $\lambda _{2}$ denotes the second-smallest of the Laplacian associated
with graph $\mathcal{G}$ and 
\[\kappa =2(L\mu_{\min}-a\lambda _{2}\mu_{\max})+\frac{4\gamma\xi}{N}(L^{2}+2a^{2}\overline{d}^{2}).\]
}

\subsection{Convergence}

We are now ready to state and prove the main theorem. 

As in \cite{ghadimi2013stochastic}, convergence is described in terms of the expected value of the average squared norm of the gradient in the first $K$-updates. The ensuing corollary goes into further detail by describing the same result in terms of {\em real-time} elapsed and not just on a total number of iterations.

\medskip \textbf{Theorem 1}: \emph{
Suppose Assumptions 0, 1, 2 and 3 hold.
Choose $\gamma<$$\min \{\bar{\gamma}_{1},\bar{\gamma}_{2}\}$, where 
\begin{align*}
\bar{\g}_1=N\frac{2a\lambda _{2}\mu_{\max}-L(2\mu_{\min}+{\xi}/{2})}{6\xi(L^{2}+2a^{2}\overline{d}^{2})} \text{ and } &\bar{\g}_2=\frac{1}{4L(2N+1)}
\end{align*}
}\emph{are positive by choosing  $ a>\frac{4\mu_{\min}L+\xi L}{4\l_2\mu_{\max}}$. With scheme \rf{flocking}, it holds that 
\begin{eqnarray*}
	&&\mathbb{E}\left[ \frac{1}{K}\sum_{k=0}^{K-1}\mathbb{E}[\Vert \tr \ell (%
	\bar{\theta}_{k})\Vert ^{2}]\right] \\
	&\leq &\frac{1}{\eta K}\left[ \ell (\bar{\theta}_{0})+L\overline{V}_{0}+%
	\frac{KL\gamma^{2}\xi\sigma ^{2}}{N^{2}}(1+\frac{1}{2N})\right],
\end{eqnarray*}
where $\eta=\frac{\gamma\xi }{N}\Big( \frac{1}{2}-2{\gamma}L(2+\frac{1}{N})\Big)$. 
}

\medskip
The regularization penalty parameter $a$
must be high enough to ensure cohesion between local models. 
This condition is weaker with a higher degree of
connectivity (i.e. higher values of $\lambda _{2}$).

Note also that for fixed $N>0$, when $a\to \infty$, then  $\g \propto 1/a$. So convergence, as characterized by Theorem 1, may be slower. This is not necessarily the case since the conditions in Theorem 1 identify a wide range of choices for $a$ and $\gamma$. For example, simulations indicate that for fixed $\gamma$ higher values of $a$ may speed up convergence (see Figure 3 (c)).

\subsection{Real-time Performance}
The analysis in Theorem 1 takes place in the time scale indexed by $k>0$ and associated with the clicks associated with a Poisson process with rate $\mu>0$. To embed the result in Theorem 1 in \emph{real-time}, recall that $%
\{D(t):t\geq 0\}$ is the counting process governing the aggregation of all
data streams. Given our assumption on computation times being negligible, the total number of
updates completed in $[0,t)$ is also $D(t)$. Let us define the  conditional average
squared gradient norm $ \Vert \bar{\tr} \ell_{t}\Vert ^{2}$ in the
interval $[0,t)$ as follows: 
\eq{ll}{
 \E[\norm{\bar{\tr} \ell_{t}} ^{2}|D(t)]\triangleq \frac{%
1}{D(t)}\sum_{k=1}^{D(t)}\Vert \tr \ell (\bar{\theta}_{k})\Vert ^{2}.}

Hence, the result in Theorem 1 can be reinterpreted { by taking expectation of \rf{ll} over $D(t)$  }as: 
\begin{align*}
\mathbb{E}[\Vert \bar{\tr} \ell _{t}\Vert ^{2}]& =\E\Big[\E[\Vert \bar{\tr} \ell_{t}\Vert ^{2}|D(t)]\Big] \\
& \leq \mathbb{E}\Big[\frac{1}{\eta D(t)}\big( \ell (\bar{\theta}_{0})+L\overline{V}_{0}\big)+%
\frac{ L\gamma^{2}\xi\sigma ^{2}}{\eta N^{2}}\big(1+\frac{1}{2N}\big)\Big] \\
& =\frac{(\ell (\bar{\theta}_{0})+L\overline{V}_{0})(1-e^{-\mu t})}{\eta \mu
t}+\frac{L{\gamma}^2\xi\sigma ^{2}}{\eta N^{2}}(1+\frac{1}{2N}).
\end{align*}%
According to Theorem 1, { and using $\g\sim \frac{1}{N}$}, the coupling of solutions via the network
regularization penalty implies the ensemble average approximates a
stationary point in the sense that: 
\begin{equation*}
\lim \sup_{t\rightarrow \infty }\mathbb{E}[\Vert \bar{\tr} \ell_{t}\Vert ^{2}]=\mathcal{O}(\frac{\sigma ^{2}}{N^{2}}).
\end{equation*}%
The \emph{approximation quality} is monotonically increasing in the number
of nodes.  The convergence properties outlined above are related to the ensemble average. It is, therefore, necessary to examine the degree to which {\em individual} models differ from the ensemble average. This is the gist of the next result.

\medskip

\textbf{Corollary 1}:
\emph{With the same assumptions and definitions in Theorem 1, it holds that 
\al{
\mathbb{E}\left[ \frac{1}{K}\sum_{k=0}^{K-1}\overline{V}_{k}\right] \leq& \frac{1}{K|\k|}\Big[\big(\frac{N}{\g}+\frac{4L\g\xi}{\eta}\big)\overline{V}_0+\frac{4\g\xi}{\eta}l(\bar{\t}_0)\Big]+
\frac{4L\g^3\xi^2\s^2}{\eta|\k|N^{2}}(1+\frac{1}{2N})+\frac{\g\xi\s^2}{|\k|N}.
}
}
{ We embed the result in Corollary 1 in \emph{real-time}. Define the conditional average of $\bar{V}_k$ in the interval $[0,t)$ as
	\[\E[\bar{V}_t|D(t)]\triangleq \frac{1}{D(t)}\sum_{k=1}^{D(t)} \bar{V}_k.   \]
The random process $\{\bar{V}_t:t>0\}$ tracks the average distance of individual models to the ensemble average.           
Similar to the discussion of Theorem 1, the \emph{real-time} result of Corollary 1 is as follows:                                                                     
\al{\E[\bar{V}_t]=& \E\Big[\E[\bar{V}_t|D(t)]\Big]\\
\leq&
 \frac{1}{D(t)|\k|}\Big[\big(\frac{N}{\g}+\frac{4L\g\xi}{\eta}\big)\overline{V}_0+\frac{4\g\xi}{\eta}l(\bar{\t}_0)\Big]+
\frac{4L\g^3\xi^2\s^2}{\eta|\k|N^{2}}(1+\frac{1}{2N})+\frac{\g\xi\s^2}{|\k|N}\\
=&\frac{1-e^{-\mu t}}{\mu t|\k|}\Big[\big(\frac{N}{\g}+\frac{4L\g\xi}{\eta}\big)\overline{V}_0+\frac{4\g\xi}{\eta}l(\bar{\t}_0)\Big]+
 +\frac{4L\g^3\xi^2\s^2}{\eta|\k|N^{2}}(1+\frac{1}{2N})+\frac{\g\xi\s^2}{|\k|N}.}
This implies the asymptotic difference between individual models and the ensemble average satisfies:
\[ \lim\sup_{t\to \infty}\E[\bar{V}_t]=\mathcal{O}(\frac{\sigma ^{2}}{N}). \]
}
{The network regularization parameter $a>0$ plays an important role in controlling the upper bound of in Corollary 1. For fixed $N>0$, when $a\to \infty$, then  $\g, \eta\propto 1/a$ and $|\k|\propto a^2$,  it follows that $\mathbb{E}\left[ \frac{1}{K}\sum_{k=0}^{K-1}\overline{V}_{k}\right]\propto 1/a$. Hence, the upper bound in Corollary 1 can be made arbitrarily small by choosing large enough $a$.}

\subsection{Comparison to Federated Learning}
\medskip 
We now present the counterpart convergence result 
regarding to FL.

\medskip
\textbf{Proposition 1}: \emph{Suppose Assumptions 0, 1, 2 and 3 hold. For scheme \ref{SG}, with a choice $\g\in (0,\frac{2}{L})$,
it holds that:
\al{\mathbb{E}\Big[\frac{1}{K}\sum_{k=0}^{K-1}\Vert \tr \ell (\theta
_{k})\Vert ^{2}\Big] 
\leq \frac{\ell (\theta _{0})}{\td{\eta} K}+\frac{L\gamma ^{2}}{2\td{\eta}}\sum_{i=1}^{N}\frac{\mu_i}{\mu}\s_i^2,}
with $\td{\eta}=$}$\gamma (1- \frac{L\gamma}{2})$%
\emph{. }

\medskip To embed the process in Proposition 1 in {\em{real-time}}, let us define
the average squared gradient norm $\Vert \tr \td{\ell}
_{t}\Vert ^{2}$ in the interval $[0,t)$ as follows: 
\eqs{p1}{ \E[\Vert \tr \td{\ell} _{t}\Vert ^{2}|D(t)]\triangleq \frac{%
		1}{D(t)}\sum_{k=1}^{D(t)}\Vert \tr \ell ( \theta_{k})\Vert ^{2}.
}
Hence, the result in Proposition 1 can be reinterpreted  by taking expectation of \rf{p1} over $D(t)$  as: 
\begin{align*}
\mathbb{E}[\Vert \tr \td{\ell} _{t}\Vert ^{2}]& =\E\Big[\E[\Vert \tr \td{\ell} _{t}\Vert ^{2}|D(t)]\Big] 
 \leq \mathbb{E}\Big[\frac{\ell (\theta _{0})}{\td{\eta} D(t)}+\frac{L{\gamma}^2}{2\td{\eta}}\sum_{i=1}^{N}\frac{\mu_i}{\mu}\s_i^2\Big] 
 =\frac{\ell (\theta_{0})(1-e^{-{\mu} t})}{\td{\eta} {\mu}
	t}+\frac{L{\gamma}^2}{2\td{\eta}}\sum_{i=1}^{N}\frac{\mu_i}{\mu}\s_i^2.
\end{align*}%
 To compare FL with NR, we also make $\g\sim \frac{1}{N}$. The asymptotic approximation quality is given by:
\begin{equation*}
\lim \sup_{t\rightarrow \infty }\mathbb{E}[\Vert\tr\td{\ell}
_{t}\Vert ^{2}]= O(\frac{1}{N}\sum_i \frac{\mu_i}{\mu} \sigma^2_i),
\end{equation*}
which suggests that the approximation quality is determined by the {\em  faster} data streams. This leads to unsatisfactory performance whenever $\mu_i \propto \sigma^2_i$ (i.e. faster data streams are also less accurate).
{Evidently, the opposite holds true when faster nodes are also more accurate, i.e. $\mu_i \propto 1/\sigma^2_i$. However, in many real-time machine learning applications, this is not likely to be the case. Obtaining higher precision gradient estimates requires larger batches and/or increased computation. Thus nodes with higher precision are less likely to be the faster ones. }

\section{Testbed: Realtime \emph{Deep} Learning}

In this section, we report the results of NR (scheme \rf{flocking})  to
distributed real-time learning from three aspects: the comparison with FL (scheme \rf{SG}),  the effects of the regularization parameter $a$, and the effects of the network connectivity.
\medskip

The specific learning task is to classify
handwritten digits between $0$ and $9$ digits as given in the MNIST data set 
\cite{lecun1998gradient}. The dataset is composed of $10000$ testing items
and $60,000$ training items. Each item in the dataset is a black-and-white
(single-channel) image of 28 $\times$ 28 pixels of a handwritten digit
between $0$ and $9$.

\medskip
  In the first two experiments, we implement schemes in a heterogeneous setting  with $5$ nodes, and  the third experiment with $20$ nodes in a  homogeneous setting.  In the test-bed MNIST streams according to independent Poisson processes. Gradient estimates are obtained with different mini-batch sizes. Evidently, a smaller mini-batch size implies noisier gradient estimates. The detailed experimental settings are summarized in Table 1. In the heterogeneous setting,  ``node 0" is the fastest and noisiest in producing gradient estimates.
 
\begin{table}[htp]
	\centering
	\begin{tabular}{c|cccc}
		\hline
		Setting & Stream ID & \# Nodes & Mini-batch Size & $\mu_i$\\
		\hline
		    & $D_0$    &1 & 1     & 8\\
	       Heterogeneous	& $D_1-D_4$&4 & 64    & 1 \\
		\hline
		 	Homogeneous & All streams& 20&4    & 1 \\
		\hline
	\end{tabular}\\

	\begin{minipage}{12cm}%
		\medskip
	\small { Table 1.The experiment hyperparameters of the two settings, including the data stream ID (Stream ID), number of nodes involved (\# Nodes), the number of images arrived as a mini-batch (Mini-batch Size), and the Poisson rate of the corresponding stream is ($\mu_i$)}.
\end{minipage}
	\label{tab:addlabel}%
\end{table}%

We use the \emph{Ray} platform (see \cite%
{moritz2018ray}) which is a popular library with shared memory supported, allowing information exchange between local nodes without copying as well as avoiding a central
bottleneck. For low-level computation, Google TensorFlow is used.
 We use a Convolutional Neural Network (CNN) with two 2D
Convolutions each with kernel size $5 \times 5$, stride 1 and 32, 64 filters. Each convolution layer is followed
by a Max-pooling with a $2\times2$ filter and stride of 2. These layers are
then followed by a Dense Layer with 256 neurons with 0.5 dropout and sigmoid activation
followed by 10 output neurons and Softmax operation. Cross entropy is used
as a performance measure (i.e. loss).\footnote{  With max-pooling the loss function is not differentiable in a set of measure zero. If in the course of execution a non-differentiable point is encountered, Tensorflow assumes a zero derivative}. Details on the implementation are available at: \url{https://github.com/wangluochao902/Network-Regularized-Approach}.

\medskip
We present the experimental results in mean plots with stand error bar. The means are computed across $10$ trials under the same hyper-parameters (namely, $\gamma$ and $a$).

\subsection{Comparison to Federated Learning}
In this experiment, we compare NR with FL in the heterogeneous setting. In Figure 2, we plot the means of the ensemble average of NR and FL with different learning rates. 
\begin{figure}[htp]
	\centering
	\subfloat[]{\label{figur:1}\includegraphics[width=80mm]{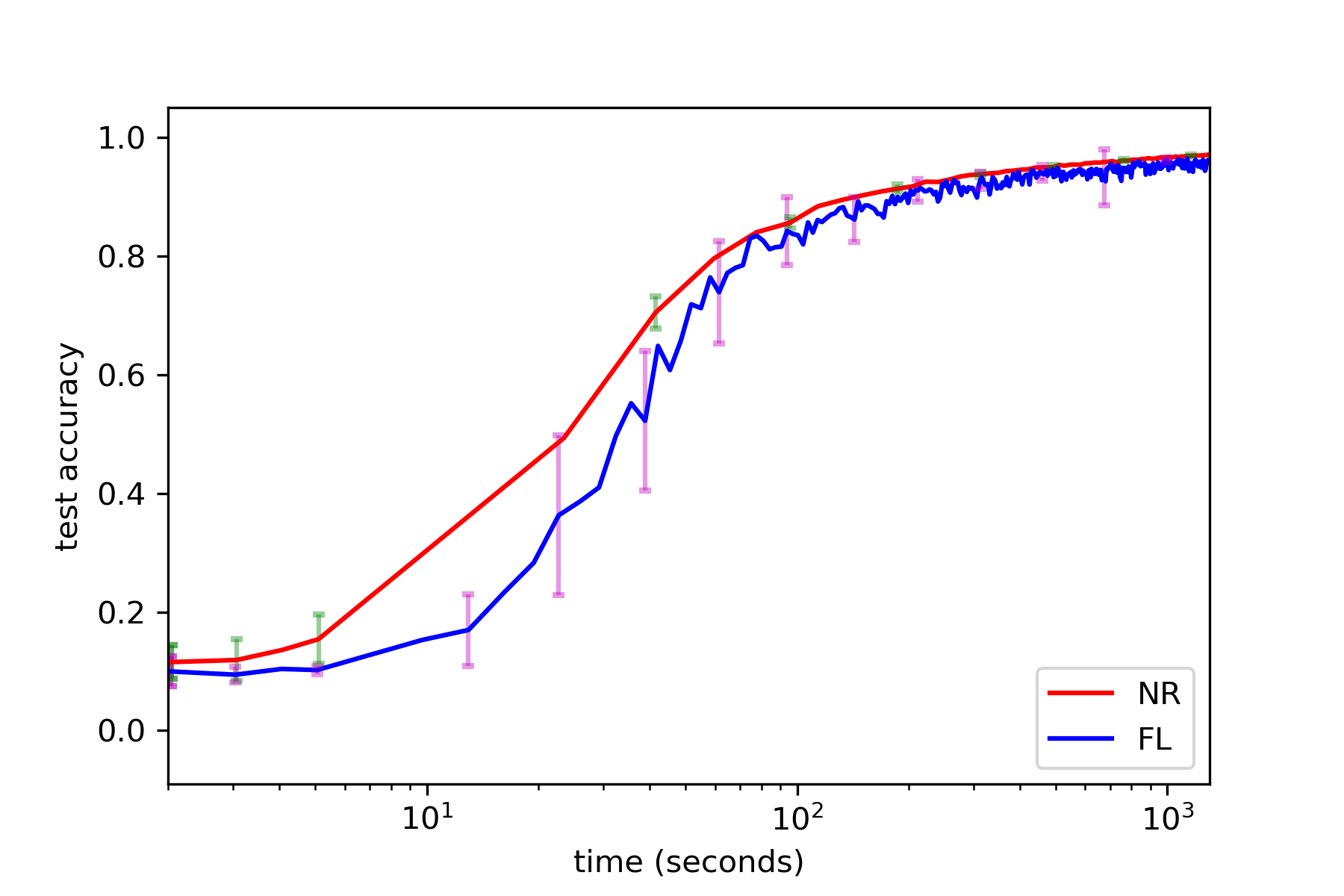}}
	\subfloat[]{\label{figur:3}\includegraphics[width=80mm]{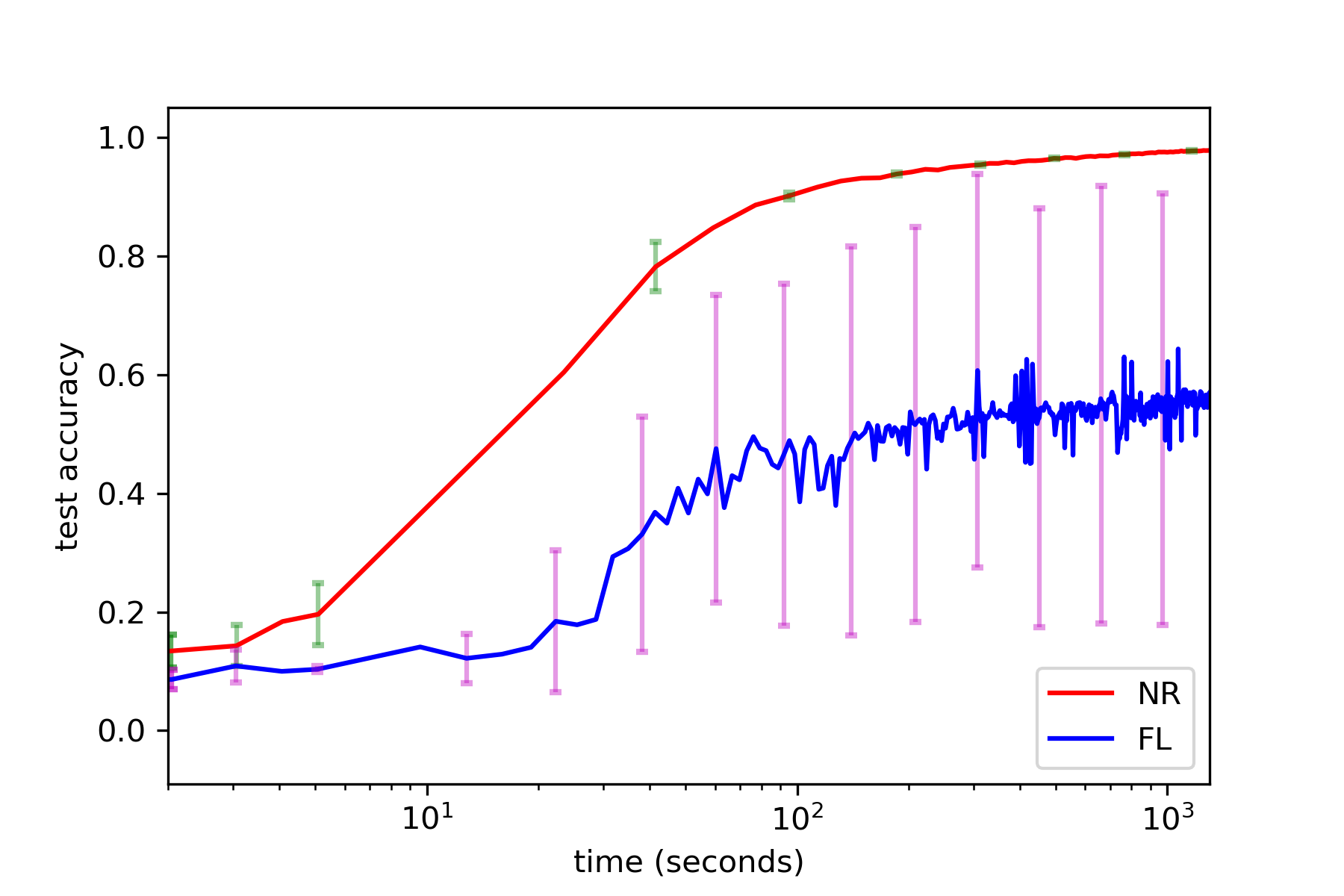}}
	\\
	\begin{minipage}{12cm}%
		\small Figure 2. The mean plot of ensemble average computed under the schemes of NR and FL in heterogeneous setting. The parameter $a$ is set to $10$ and  the network is fully connected. The learning rate $\g$ is set  to $0.002$ in (a) and $0.004$ in (b). 
	\end{minipage}
\end{figure}

We can observe from Figure (a) that when the learning rate is moderate, both FL and  NR  can converge, but the empirical standard deviation of FL is much larger than that of NR. With increased $\g$,  FL fails to converge while NR still performs relatively well, as shown in Figure (b). We can see that NR is more robust with respect to the learning rate.

\medskip
\subsection{The Effects of Regularization Parameter}

In this experiment, we look at the effects of changing the regularization parameter $a$. In Figure 3, we present the means of each node as well as the ensemble average.
\begin{figure}[htp]
	\centering
	\subfloat[]{\label{figur:1}\includegraphics[width=80mm]{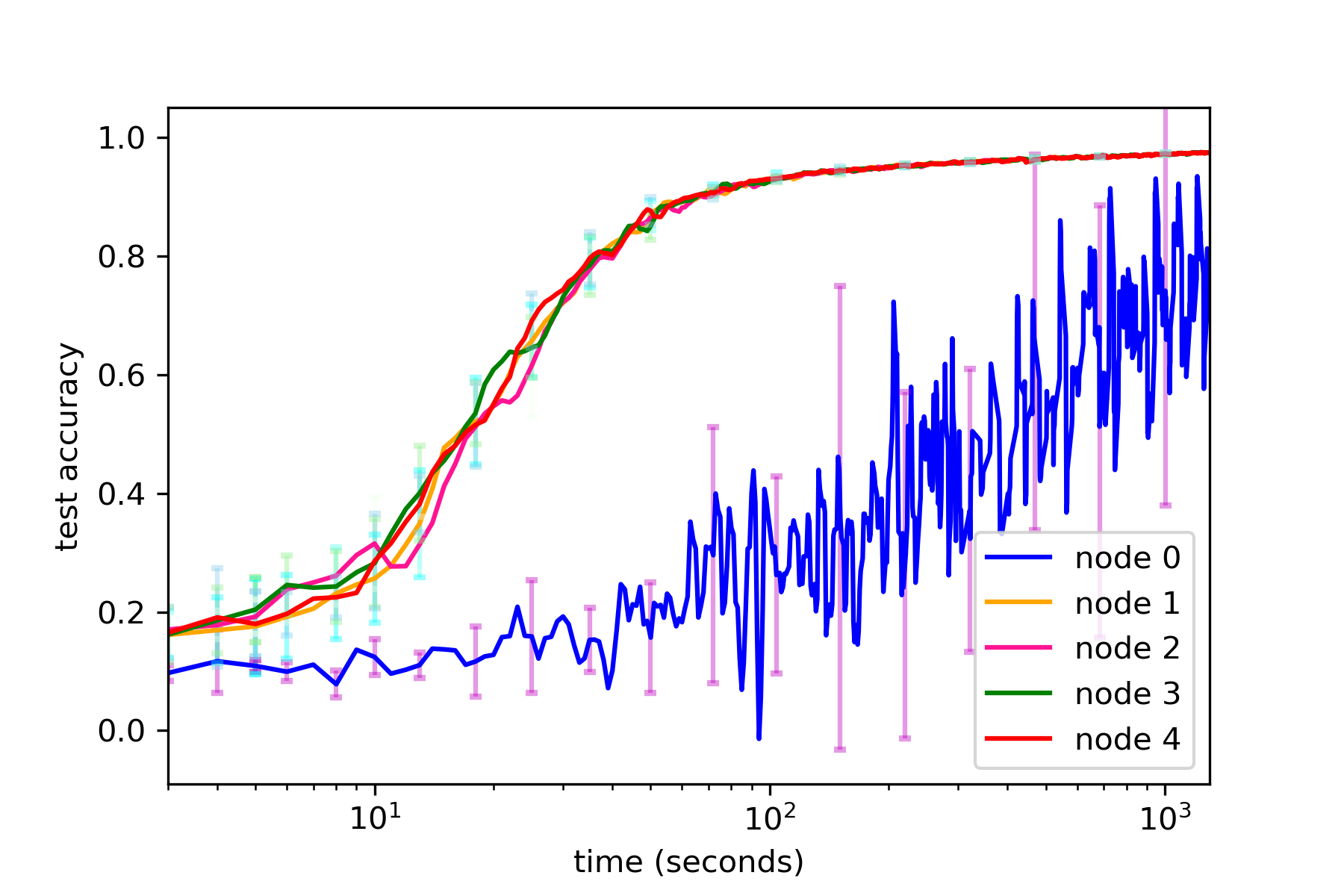}}
	\subfloat[]{\label{figur:3}\includegraphics[width=80mm]{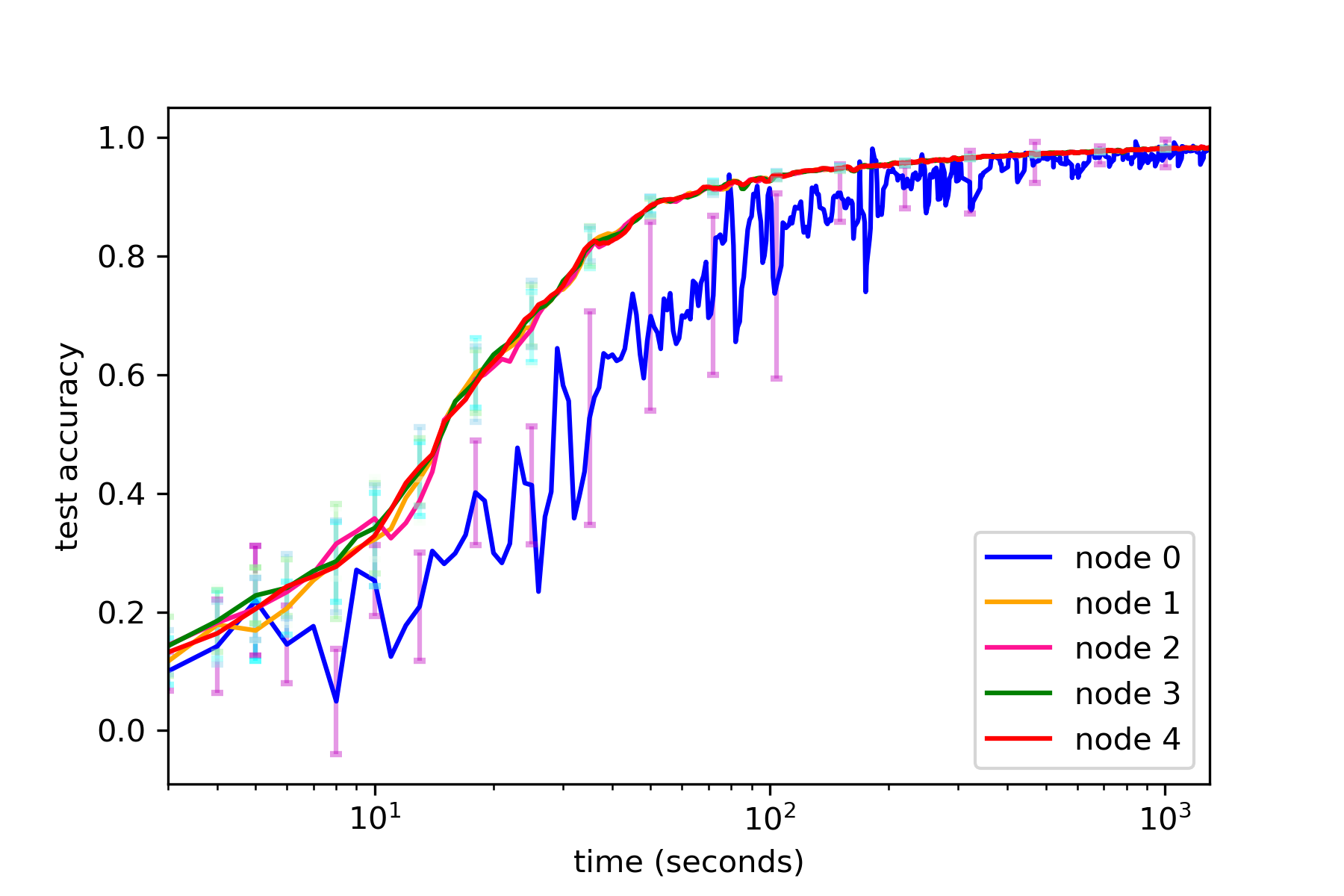}}\\
	\subfloat[]{\label{figur:3}\includegraphics[width=80mm]{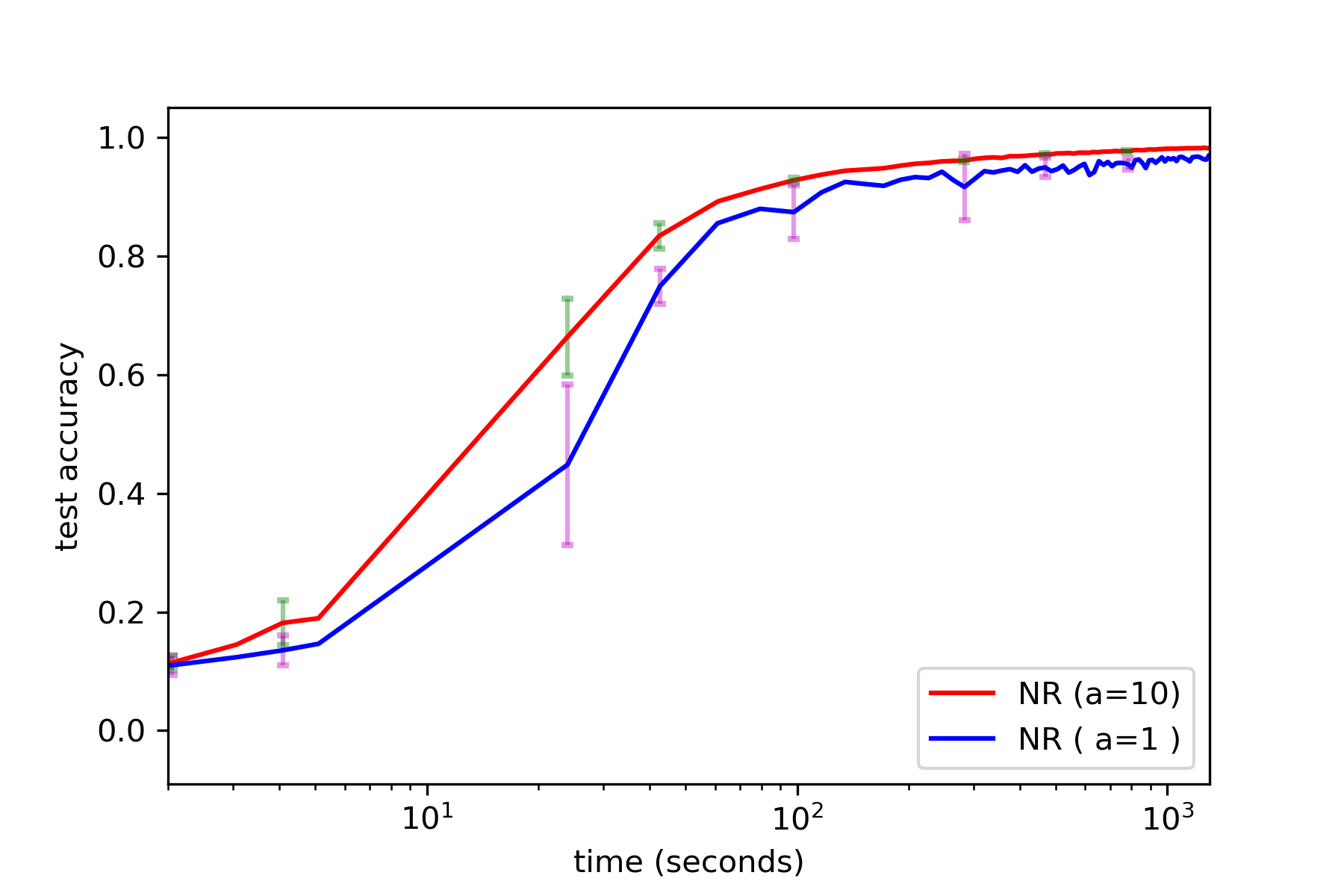}}
	\\
	\begin{minipage}{12cm}%
		\small Figure 3. The mean plot of each node computed under the scheme of NR in heterogeneous setting. The parameter $\g$ is set to $0.01$ and the network is fully connected.  The regularization parameter $a$ is set to $1$ in (a) and $10$ in (b). 
	  The mean plot of the  ensemble average under two choices of $a$ is presented in (c).
	\end{minipage}
\end{figure}

As we increase $a$ from $1$ to $10$, we can observe from Figure 3 (a) and (b) that the consensus among nodes increases and the empirical mean standard deviation of the ``node 0" decreases. As presented  in Corollary 1, the regularization parameter $a$ influences the degree of similarity between individual models and the ensemble average.  Note that we only identify a range of values for  $a$ (lower bound) and $\g$ (upper bound) for which convergence is guaranteed so that  a higher value of $a$ does not necessarily imply slower convergence, as shown in Figure 3 (c).

\subsection{The Effects of Network Connectivity}

In the third experiment, we check the effect of increased connectivity in  the homogeneous setting by
using a Watts-Strogatz ``small world" topology (see \cite{Watts}), in which each node is connected with 2 (or 8) nearest
neighbors.

\begin{figure}[htp]
	\centering
	\subfloat[]{\label{figur:1}\includegraphics[width=80mm]{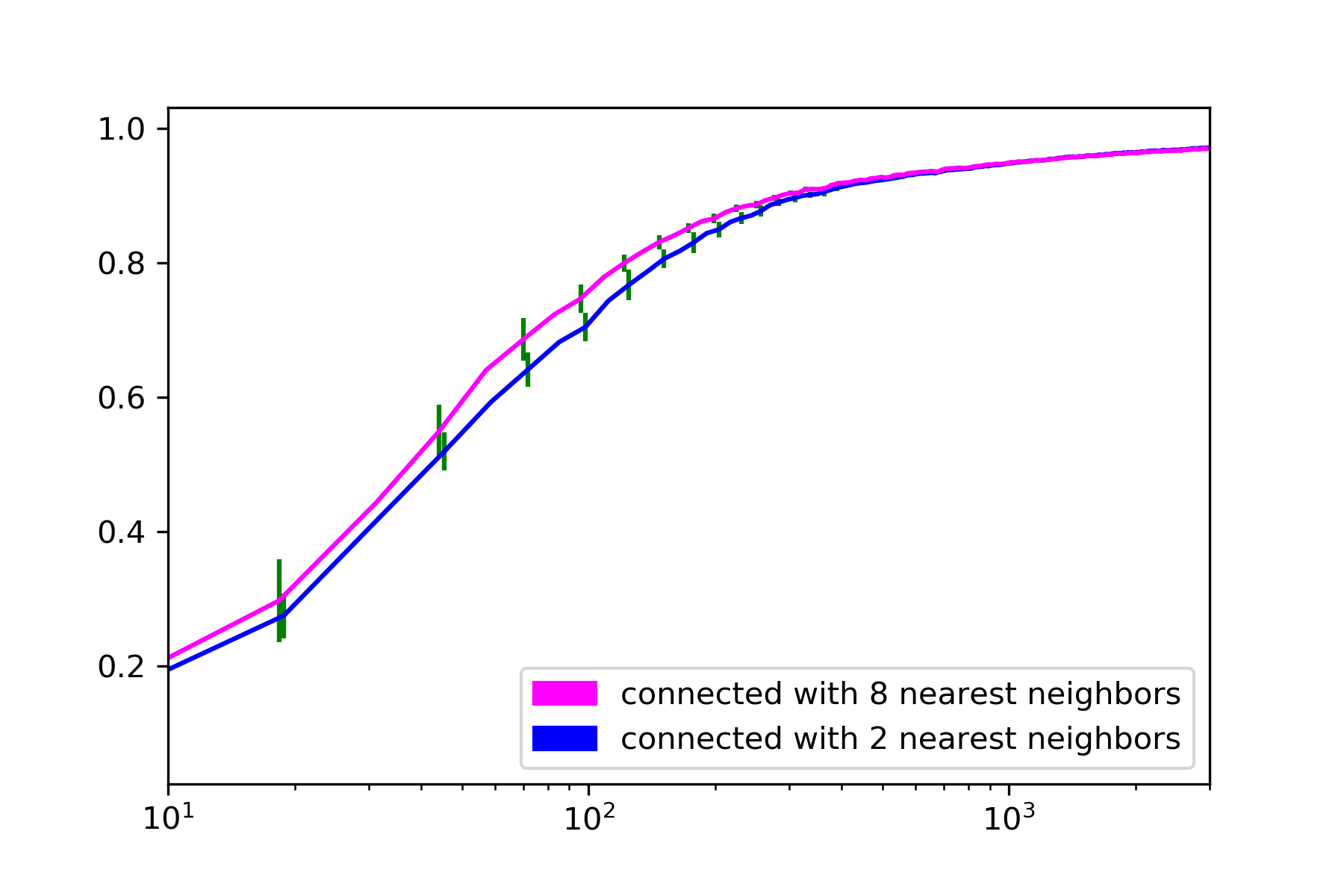}}\\

	\begin{minipage}{8cm}%
		\small Figure 4. The mean plot of ensemble average computed under the scheme of NR in the homogeneous setting. The parameter $a$ is set to $10$ and the learning rate $\g$ is set to  $0.001$. 
	\end{minipage}
\end{figure}

We can see from Figure 4 that increasing the connectivity of the topology only improves the performance slightly, meaning that only a limited connectivity is needed for the network regularized approach to enjoy a satisfactory rate of convergence.

\section{Conclusions}

In many application domains, data streams through a network of heterogeneous nodes in different geographic locations. When there is high data payload (e.g. high-resolution video), assembling a diverse batch of data
points in a central processing location in order to update a model entails
significant latency. In such cases, a distributed architecture for learning relying on a network of interconnected ``local" nodes may prove advantageous.
We have analyzed a distributed scheme in which every local node implements stochastic gradient updates every time a data
point is obtained. To ensure robust estimation, a local regularization penalty is
used to maintain a measure of cohesion in the ensemble of models. We show
the ensemble average approximates a stationary point.
The approximation quality is superior to that of FL, especially when there is heterogeneity in gradient estimation quality. We also show that our approach is robust against changes in the learning rate and network connectivity.  
We illustrate the
results with an application to deep learning with convolutional neural
networks.

In future work we plan to study different localized model averaging schemes. A careful selection of weights for computing local average model ensures a reduction of estimation variance. This is motivated by the literature on the optimal combination of forecasts (see \cite{Bates}). For example, weights minimizing the sample mean square prediction error are of the form  $\frac{\hat{\sigma}^{-2}_i}{\sum_{j=1}^{N} \hat{\sigma}^{-2}_j}$ where $\hat{\sigma}^{2}_i$ is the estimated mean squared prediction error of the $i$-th model.

\section{Appendix}

\textbf{Proof of Lemma 1}

Note that%
\begin{align*}
\bar{\theta}_{k+1}& =\frac{1}{N}\sum_{i=1}^{N}\left[ \theta _{i,k}-\gamma\mathbf{1}_{i,k}[\tr \ell (\theta _{i,k})+a\tr \mathcal{F}%
_{i,k}+\varepsilon _{i,k}]\right] \\
& =\bar{\theta}_{k}-\frac{\g}{N}\sum_{i=1}^{N}\tr \ell (\theta
_{i,k})\mathbf{1}_{i,k}  -\frac{a\g}{N}\sum_{i=1}^{N}\tr \mathcal{F}_{i,k}\mathbf{1}%
_{i,k}-\frac{\g}{N}\sum_{i=1}^{N}\varepsilon _{i,k}\mathbf{1}_{i,k}.
\end{align*}%
Hence, $e_{i,k+1}=\theta _{i,k+1}-\bar{\theta}_{k+1}%
=e_{i,k}-\gamma \delta _{i,k}$. Then%
\al{V_{i,k+1}& =(e_{i,k}-\gamma \d_{i,k})\trs(e_{i,k}-\gamma\delta_{i,k}) \\
& =e_{i,k}\trs e_{i,k}-2\gamma e_{i,k}\trs\delta _{i,k}+\gamma^{2}\norm{ \d_{i,k}}^{2} \\
& =V_{i,k}-2\gamma e_{i,k}\trs(\d_{i,k}^{f}+\d_{i,k}^{g}+\d_{i,k}^{n})+\g^2\norm{ \d_{i,k} }^{2},}
and%
\begin{equation*}
\overline{V}_{k+1}=\overline{V}_{k}-\frac{2\gamma}{N}\sum_{i=1}^{N}e_{i,k}\trs(\delta _{i,k}^{f}+\delta _{i,k}^{g}+\delta _{i,k}^{n})+\frac{%
\gamma ^2}{N}\sum_{i=1}^{N}\Vert \delta _{i,k}\Vert ^{2}.
\end{equation*}%
Finally, note that%
\begin{align*}
\sum_{i=1}^{N}e_{i,k}\trs\delta _{i,k}^{f}& =\sum_{i=1}^{N}e_{i,k}\trs\left[
\tr \ell (\theta _{i,k})\mathbf{1}_{i,k}-\tr\bar{\ell} _{k}\right] 
=\sum_{i=1}^{N}e_{i,k}\trs\tr \ell (\theta _{i,k})\mathbf{1}_{i,k}, \\
\sum_{i=1}^{N}e_{i,k}\trs\delta _{i,k}^{g}&
=a\sum_{i=1}^{N}e_{i,k}\trs(\tr \mathcal{F}_{i,k}\mathbf{1}_{i,k}-%
\tr\bar{\@F}_{k}) 
=a\sum_{i=1}^{N}e_{i,k}\trs\tr \mathcal{F}_{i,k}\mathbf{1}_{i,k}, \\
\sum_{i=1}^{N}e_{i,k}\trs\delta _{i,k}^{n}&
=\sum_{i=1}^{N}e_{i,k}\trs(\varepsilon _{i,k}\mathbf{1}_{i,k}-\frac{1}{N}%
\sum_{j=1}^{N}\varepsilon _{j,k}\mathbf{1}_{j,k}) 
 =\sum_{i=1}^{N}e_{i,k}\trs\varepsilon _{i,k}\mathbf{1}_{i,k}.
\end{align*}%
So the result follows. $\blacksquare $

\medskip \textbf{Proof of Lemma 2}

In light of Lemma 1 we have:%
\begin{eqnarray*}
\mathbb{E}[\overline{V}_{k+1}|\*{\t}_k] &=&\overline{V}_{k}-\frac{2\gamma 
}{N}\sum_{i=1}^{N}\frac{\mu_i}{\mu}e_{i,k}\trs\left[ \tr \ell (\theta _{i,k})+a\tr 
\mathcal{F}_{i,k}\right] +\frac{\g^2}{N}\sum_{i=1}^{N}\Vert \delta _{i,k}\Vert ^{2}.
\end{eqnarray*}%
Let $\mathbf{e}_{k}=\left[ 
\begin{array}{cccc}
e_{1,k}\trs,  e_{2,k}\trs, \cdots, e_{N,k}\trs%
\end{array}%
\right] \trs$ and $\mathcal{L}=[l_{ij}]$ be the Laplacian matrix associated
with the adjacency matrix $A$, where $l_{ii}=\sum_{j}a_{ij}$ and $%
l_{ij}=-a_{ij}$ when $i\neq j$. For an undirected graph, the Laplacian
matrix is symmetric positive semi-definite. It follows that%
\begin{align*}
\sum_{i=1}^{N}e_{i,k}\trs\tr \mathcal{F}_{i,k}&
=\sum_{i=1}^{N}\sum_{j=1,j\neq i}^{N}\alpha _{ij}e_{i,k}\trs(e_{i,k}-e_{j,k})
 =-\sum_{i=1}^{N}\sum_{j\neq i}^{N}l_{ij}e_{i,k}\trs(e_{i,k}-e_{j,k}) \\
& =\sum_{i=1}^{N}\sum_{j\neq i}^{N}l_{ij}e_{i,k}\trs e_{j,k} 
 =\mathbf{e}_{k}\trs(\mathcal{L}\otimes I_{p})\mathbf{e}_{k}  \geq \lambda _{2}\sum_{i=1}^{N}\left\Vert e_{i,k}\right\Vert ^{2},
\end{align*}

where $\lambda _{2}:=\lambda _{2}(\mathcal{L})$ is the second-smallest
eigenvalue of $\mathcal{L}$, also called the \emph{algebraic connectivity}
of $\mathcal{G}$ \cite{godsil2013algebraic}. Thus,%
\eqs{Lemma2.1}{
\mathbb{E}[\overline{V}_{k+1}|\*{\t}_k] \leq& \overline{V}_{k}-\frac{
2\gamma }{N}\sum_{i=1}^{N}\frac{\mu_i}{\mu}e_{i,k}\trs\tr \ell (\theta _{i,k})-
\frac{2a\lambda _{2}\gamma }{N}\sum_{i=1}^{N}\frac{\mu_i}{\mu}\left\Vert e_{i,k}\right\Vert ^{2}
 +\frac{\gamma^{2}}{N}\sum_{i=1}^{N}\mathbb{E}[\Vert \delta
_{i,k}\Vert ^{2}|\*{\t}_k]  \notag \\
 =& \overline{V}_{k}-\frac{2\gamma }{N}\sum_{i=1}^{N}\frac{\mu_i}{\mu}(\tr \ell (\theta
_{i,k})-\tr \ell (\bar{\theta}_{k}))\trs e_{i,k}  -\frac{2a\lambda _{2}\gamma }{N}\sum_{i=1}^{N}\frac{\mu_i}{\mu}\left\Vert
e_{i,k}\right\Vert ^{2}+\frac{\gamma^{2}}{N}\sum_{i=1}^{N}\mathbb{E}%
[\Vert \delta _{i,k}\Vert ^{2}|\*{\t}_k]. }
By Cauchy–Schwarz inequality and Assumption $3$, we can obtain that
\al{-(\tr \ell (\theta_{i,k})-\tr \ell (\bar{\theta}_{k}))\trs e_{i,k}&\leq \norm{\tr \ell (\theta_{i,k})-\tr \ell (\bar{\theta}_{k})}\norm{e_{i,k}}\leq L \left\Vert e_{i,k}\right\Vert ^{2}.}
Define $\bar{\mu}=\mu/N$, and by the inequalities  $ \frac{\mu_{\min}}{N \bar{\mu}} \leq\frac{\mu_i}{\mu} \leq \frac{\mu_{\max}}{N \bar{\mu}}$, we can obtain
\eqs{Lemma2.1}{
\mathbb{E}[\overline{V}_{k+1}|\*{\t}_k]\leq& (1+\frac{2\g}{N\bar{\mu}}(L\mu_{\max}-a\lambda_{2}\mu_{min}))\overline{V}_{k} +\frac{\gamma^{2}}{N}\sum_{i=1}^{N}\mathbb{E}[\norm{\delta _{i,k}}^{2}|\*{\t
}_{k}].}  
We now simplify the last term in the right hand side of \rf{Lemma2.1}.
First we note that: 
\begin{equation}
\mathbb{E}[\Vert \delta _{i,k}\Vert ^{2}|\*{\t}_k]=\mathbb{E}[\Vert
\delta _{i,k}^{f}+\delta _{i,k}^{g}\Vert ^{2}|\*{\t}_k]+\mathbb{E}[\Vert
\delta _{i,k}^{n}\Vert ^{2}|\*{\t}_k].  \label{Lemma2.2}
\end{equation}%
The first term in the right hand side of (\ref{Lemma2.2}) can be further
described as follows:%
\al{& \gamma^{2}\mathbb{E}[\Vert \delta _{i,k}^{f}+\delta _{i,k}^{g}\Vert
^{2}|\*{\t}_k] \\
=& \gamma^{2}\mathbb{E}\Big[\left. \left\Vert \tr \ell (\theta _{i,k})%
\mathbf{1}_{i,k}-\tr\bar{\ell} _{k}+a(\tr \mathcal{F}_{i,k}\mathbf{1}%
_{i,k}-\tr\bar{\@F}_{k})\right\Vert ^{2}\right\vert \*{\t}_k\Big]
\\
=& \gamma^{2}\mathbb{E}\Big[\Vert (1-\frac{1}{N})[\tr \ell (\theta
_{i,k})+a\tr \mathcal{F}_{i,k}]\mathbf{1}_{i,k}+  \frac{1}{N}\sum_{j\neq i}^{N}[\tr \ell (\theta
_{j,k})+a\tr \mathcal{F}_{j,k}]\mathbf{1}_{j,k}\Vert ^{2}|\*{\t}_k\Big] \\
=& \frac{\gamma ^{2}}{N}\Big[(1-\frac{1}{N})^{2}\frac{\mu_i}{\bar{\mu}}%
\left\Vert \tr \ell (\theta _{i,k})+a\tr \mathcal{F}_{i,k}\right\Vert
^{2}+ \frac{1}{N^2}\sum_{j\neq i}^{N}\frac{\mu _{j}}{\bar{\mu}}\left\Vert \tr \ell (\theta _{j,k})\mathbf{+%
}a\tr \mathcal{F}_{j,k}\right\Vert ^{2}\Big].}
This leads to: 
\begin{align}
 \sum_{i=1}^{N}\gamma^{2}\mathbb{E}[\Vert \delta _{i,k}^{f}+\delta
_{i,k}^{g}\Vert ^{2}|\*{\t}_k]  &  \leq \frac{\gamma ^{2}\xi }{N}\Big[(1-\frac{1}{N}%
)^{2}\sum_{i=1}^{N}\left\Vert \tr \ell (\theta _{i,k})\mathbf{+}a\tr 
\mathcal{F}_{i,k}\right\Vert ^{2} \notag  +\frac{1}{N^2}\sum_{i=1}^{N}\sum_{j\neq i}^{N}\left\Vert \tr \ell (\theta _{j,k})\mathbf{+}a\tr 
\mathcal{F}_{j,k}\right\Vert ^{2}\Big]  \notag \\
& \leq \frac{\g^2\xi}{N}\sum_{i=1}^{N}\left\Vert \tr \ell (\theta _{i,k})\mathbf{+}a\tr 
\mathcal{F}_{i,k}\right\Vert ^{2}.  \label{Lemma2.3}
\end{align}%
Finally,%
\begin{align}
 \gamma^{2}\sum_{i=1}^{N}\mathbb{E}[\Vert \delta _{i,k}^{n}\Vert
^{2}|\*{\t}_k]  
& =\gamma^{2}\sum_{i=1}^{N}\mathbb{E}\Big[\Vert (1-\frac{1}{N})\varepsilon
_{i,k}\mathbf{1}_{i,k}-\frac{1}{N}\sum_{j\neq i}^{N}\varepsilon _{j,k}%
\mathbf{1}_{j,k}\Vert ^{2}|\*{\t}_k\Big] \notag \\
& =\frac{\gamma ^{2}}{N}\Big[(1-\frac{1}{N})^{2}\sum_{i=1}^{N}\frac{\mu_i }{\bar{\mu}}\mathbb{E}[\Vert \varepsilon _{i,k}\Vert ^{2}|\*{\t}_k]  +\frac{1}{N^{2}}\sum_{i=1}^{N}\sum_{j\neq i}^{N}\frac{%
\mu_i }{\bar{\mu}}\mathbb{E}[\Vert \varepsilon
_{j,k}\Vert ^{2}|\*{\t}_k]\Big]  \notag \\
&  \leq \frac{\gamma ^{2}}{N}\big(\frac{\mu _{\max }}{%
\mu _{\min }}\big)\sum_{i=1}^{N}\mathbb{E}[\Vert \varepsilon _{i,k}\Vert ^{2}|\*{\t}_{k}]\leq \frac{  \g^2\xi\sigma ^{2}%
}{N}.  \label{Lemma2.4}
\end{align}%
We use inequalities (\ref{Lemma2.3}) and (\ref{Lemma2.4}) with (\ref%
{Lemma2.2}) to obtain an upper bound of (\ref{Lemma2.1}) as follows: 
\eqs{Lemma2.5}{\mathbb{E}[\overline{V}_{k+1}|\*{\t}_k]  
\leq &(1+\frac{2\gamma }{N}(L\mu_{\min}-a\lambda _{2}\mu_{\max}))\overline{V}_{k} +\frac{\gamma^{2}\xi}{N^{2}}\sum_{i=1}^{N}\Vert \tr \ell (\theta
_{i,k})+a\tr \mathcal{F}_{i,k}\Vert ^{2}+\frac{\gamma^{2}\xi\sigma^{2}}{N^{2}}.  }
Finally, we
analyze the third term on the right hand side of (\ref{Lemma2.5}). By
Parallellogram law 
\begin{align*}
 \Vert \tr \ell (\theta _{i,k})+a\tr \mathcal{F}_{i,k}\Vert ^{2} &=2\Vert \tr \ell (\theta _{i,k})\Vert ^{2}+2\Vert a\tr \mathcal{F}%
_{i,k}\Vert ^{2}-\Vert \tr \ell (\theta _{i,k})-a\tr \mathcal{F}%
_{i,k}\Vert ^{2} \\
& \leq 2\Vert \tr \ell (\theta _{i,k})\Vert ^{2}+2\Vert a\tr \mathcal{F%
}_{i,k}\Vert ^{2}
\end{align*}%
In addition, 
\begin{align*}
\Vert \tr \mathcal{F}_{i,k}\Vert ^{2}& =\deg (i)^{2}\left\Vert
\sum_{j=1,j\neq i}^{N}\frac{\alpha _{ij}(\theta _{i,k}-\theta _{j,t})}{\deg
(i)}\right\Vert ^{2} \\
& \leq \deg (i)\sum_{j=1,j\neq i}^{N}\alpha _{ij}\left\Vert \theta
_{i,k}-\theta _{j,k}\right\Vert ^{2}  \leq \bar{d}\sum_{j=1,j\neq i}^{N}\alpha _{ij}\left\Vert \theta
_{i,k}-\theta _{j,k}\right\Vert ^{2}.
\end{align*}%
which implies 
\begin{align*}
\sum_{i=1}^{N}\Vert \tr \mathcal{F}_{i,k}\Vert ^{2}& \leq \bar{d}%
\sum_{i=1}^{N}\sum_{j=1,j\neq i}^{N}\alpha _{ij}\left\Vert \theta
_{i,k}-\theta _{j,k}\right\Vert ^{2}  =\bar{d}\sum_{i=1}^{N}\sum_{j\neq i}^{N}\alpha _{ij}\left\Vert
e_{i,k}-e_{j,k}\right\Vert ^{2} \\
& \leq 2\bar{d}\sum_{i=1}^{N}\sum_{j\neq i}^{N}\alpha _{ij}(\left\Vert
e_{i,k}\right\Vert ^{2}+\left\Vert e_{j,k}\right\Vert ^{2})  \leq 4\bar{d}^{2}\sum_{i=1}^{N}\left\Vert e_{i,k}\right\Vert ^{2}=4N\bar{d}%
^{2}\overline{V}_{k}.
\end{align*}%
Thus,%
\eqs{Lemma2.6}{
\sum_{i=1}^{N}\Vert \tr \ell (\theta _{i,k})+a\tr \mathcal{F}%
_{i,k}\Vert ^{2}  
& \leq 2\sum_{i=1}^{N}\Vert \tr \ell (\theta _{i,k})\Vert ^{2}+8a^{2}N%
\overline{d}^{2}\overline{V}_{k} \\
& =4N\left\Vert \tr \ell (\bar{\theta}_{k})\right\Vert
^{2}+4\sum_{i=1}^{N}\left\Vert \tr \ell (\theta _{i,k})-\tr \ell (\bar{%
\theta}_{t})\right\Vert ^{2}+8a^{2}N\overline{d}^{2}\overline{V}_{k}  \\
& \leq 4N\left\Vert \tr \ell (\bar{\theta}_{k})\right\Vert
^{2}+4N(L^{2}+2a^{2}\overline{d}^{2})\overline{V}_{k}.}
The result follows by using the previous inequality to obtain an upper bound
for the right hand side of (\ref{Lemma2.5}). $\blacksquare $

\medskip \textbf{Proof of Theorem 1}

By Taylor expansion and Lipschitz assumption: 
\begin{align*}
\ell (\bar{\theta}_{k+1}) \leq& \ell (\bar{\theta}_{k})+\tr \ell (\bar{%
\theta}_{k})\trs(\bar{\theta}_{k+1}-\bar{\theta}_{k})+\frac{L}{2}\left\Vert 
\bar{\theta}_{k+1}-\bar{\theta}_{k}\right\Vert ^{2} \\
 =&\ell (\bar{\theta}_{k})-\frac{\g}{N}\sum_{i=1}^{N}\tr \ell (%
\bar{\theta}_{k})\trs\tr \ell (\theta _{i,k})\mathbf{1}_{i,k}  -\frac{a\g}{N}\sum_{i=1}^{N}\tr \ell (\bar{\theta}
_{k})^{T}\tr \mathcal{F}_{i,k}\mathbf{1}_{i,k} \\&-\frac{\g}{N}\sum_{i=1}^{N}\tr \ell (\bar{\theta}%
_{k})^{T}\varepsilon _{i,k}\mathbf{1}_{i,k}+\frac{L}{2}\left\Vert \bar{\theta%
}_{k+1}-\bar{\theta}_{k}\right\Vert ^{2}.
\end{align*}%

Since $\sum_{i=1}^{N}\tr \mathcal{F}_{i,k}=0$, it follows that 

\begin{align}
 E[\left. \ell (\bar{\theta}_{k+1})\right\vert \*{\t}_k] 
\leq& \ell (\bar{\theta}_{k})-\frac{\gamma \xi}{N^{2}}\sum_{i=1}^{N}\tr
\ell (\bar{\theta}_{k})^{T}\tr \ell (\theta _{i,k})+\frac{L}{2}E[
\left\Vert \bar{\theta}_{k+1}-\bar{\theta}_{k}\right\Vert ^{2}|
\*{\t}_k] \label{Theorem1.1}  \\
 \leq& \ell (\bar{\theta}_{k})-\frac{\gamma\xi }{N^{2}}\sum_{i=1}^{N}\tr
\ell (\bar{\theta}_{k})^{T}[\tr \ell (\theta _{i,k})-\tr \ell (\bar{%
\theta}_{k})] \notag  -\frac{\gamma\xi }{N}\left\Vert \tr \ell (%
\bar{\theta}_{k})\right\Vert ^{2}+\frac{L}{2}E[\left\Vert \bar{\theta}_{k+1}-\bar{\theta}%
_{k}\right\Vert ^{2}| \*{\t}_k]. 
\end{align}

Using \rf{Lemma2.6} from the proof of Lemma 2 we obtain%
\begin{align}
 E[\left. \left\Vert \bar{\theta}_{k+1}-\bar{\theta}_{k}\right\Vert
^{2}\right\vert \*{\t}_k]  
 =&\frac{\g^2}{N^{3}}\sum_{i=1}^{N}\frac{\mu _{i}}{\bar{\mu} }%
\left\Vert \tr \ell (\theta _{i,k})+a\tr \mathcal{F}_{i,k}\right\Vert
^{2} \notag 
+\frac{\g^2}{N^{3}}\sum_{i=1}^{N}\frac{\mu _{i}}{\bar{\mu} }\mathbb{E%
}[\Vert \varepsilon _{i,k}\Vert ^{2}|\*{\t}_k] \\
 \leq &\frac{4\gamma ^{2}\xi}{N^{2}}
\Big[\left\Vert \tr \ell (\bar{\theta}_{k})\right\Vert ^{2}+(L^{2}+2a^{2}%
\overline{d}^{2})\overline{V}_{k}\Big]+\frac{\gamma ^{2}\xi\sigma ^{2}}{N^{3}}. 
\label{Theorem1.2}
\end{align}
Also, 

\begin{align}
 -\tr \ell (\bar{\theta}_{k})^{T}[\tr \ell (\theta _{i,k})-\tr
\ell (\bar{\theta}_{k})]  \notag 
=& \frac{1}{2}\left\Vert \tr \ell (\bar{\theta}_{k})\right\Vert ^{2}+%
\frac{1}{2}\left\Vert \tr \ell (\theta _{i,k})-\tr \ell (\bar{\theta}%
_{k})\right\Vert ^{2}-\left\Vert \tr \ell (\theta _{i,k})\right\Vert ^{2} \notag\\
 \leq &\frac{1}{2}\left\Vert \tr \ell (\bar{\theta}_{k})\right\Vert ^{2}+%
\frac{L^{2}}{2}\left\Vert \theta _{i,k}-\bar{\theta}_{k}\right\Vert ^{2}.
\label{Theorem1.3}
\end{align}

Substituting (\ref{Theorem1.3}) and (\ref{Theorem1.2}) into (\ref{Theorem1.1}%
) we obtain: 
\begin{align}
 E[\left. \ell (\bar{\theta}_{k+1})\right\vert \*{\t}_k]  \notag 
\leq& \ell (\bar{\theta}_{k})-\frac{\gamma\xi }{2N}\left\Vert \tr \ell (%
\bar{\theta}_{k})\right\Vert ^{2}+\frac{L^{2}\gamma\xi }{2N}\bar{V}_{k} \notag \\
& +\frac{2L\gamma ^{2}\xi}{N^{2}}\Big[\left\Vert
\tr \ell (\bar{\theta}_{k})\right\Vert ^{2}+(L^{2}+2a^{2}\overline{d}^{2})%
\overline{V}_{k}\Big]+\frac{L\gamma ^{2}\xi\sigma ^{2}}{2N^{3}}.
\label{Theorem1.4}
\end{align}

Consider the function $\ell (\bar{\theta}_{k})+L\overline{V}_{k}$. From the
inequalities in (\ref{Theorem1.4}) and Lemma 2 we obtain:
\begin{equation*}
\mathbb{E}[\overline{V}_{k+1}|\*{\t}_k]\leq (1+\frac{\kappa \gamma }{N})%
\overline{V}_{k}+\frac{4\gamma^{2}\xi}{N}\left\Vert \tr \ell (\bar{%
\theta}_{k})\right\Vert ^{2}+\frac{\gamma^{2}\xi\sigma ^{2}}{N^{2}},
\end{equation*}%
\begin{align*}
 \mathbb{E}[\ell (\bar{\theta}_{k+1})+L\overline{V}_{k+1}|\*{\t}_{k}] 
 \leq& (\ell (\bar{\theta}_{k})+L\overline{V}_{k}) -\frac{\gamma\xi }{N}\left( \frac{1}{2}-2\gamma L(2+\frac{1}{N})\right) \Vert \tr \ell (\bar{\theta}_{k})\Vert ^{2} \\
& +\Big[\kappa +\frac{L\xi}{2}+\frac{2\g\xi}{N}(L^{2}+2a^{2}\overline{d}^{2})\Big]\frac{L\gamma }{N}\overline{V}_{k} +\frac{L\g^{2}\xi\sigma ^{2}}{N^{2}}(1+\frac{1}{2N}).
\end{align*}

By choosing $a>\frac{4\mu_{\min}L+\xi L}{4\l_2\mu_{\max}}$, $\bar{\g}_1>0$.
Given the choice $\g<\bar{\gamma}_{1}$ in the statement of Theorem
1, we have 
\al{\kappa +\frac{L\xi}{2}+\frac{2\g\xi}{N}(L^{2}+2a^{2}\overline{d}^{2})
=&-2a\lambda _{2}\mu_{\max}+L(2\mu_{\min}+\frac{\xi}{2})+ \frac{6\g\xi}{N}(L^{2}+2a^{2}\overline{d}^{2})\leq 0}

It follows that 
\begin{align*}
\frac{\gamma\xi }{N}\Big( \frac{1}{2}-2{\gamma}L(2+\frac{1}{N})\Big) \Vert \tr \ell (\bar{\theta}_{k})\Vert ^{2} 
\leq \ell (\bar{\theta}_{k})+L\overline{V}_{k}-\mathbb{E}\left[ \ell (\bar{%
\theta}_{k+1})+L\overline{V}_{k+1}|\*{\t}_k\right] +\frac{L{\gamma}^{2}\xi\sigma ^{2}}{N^{2}}(1+\frac{1}{2N}).
\end{align*}%
Let $\eta=\frac{\gamma\xi }{N}\Big( \frac{1}{2}-2{\gamma}L(2+\frac{1}{N})\Big)$.
By definition ${\gamma}<\bar{\gamma}_{2}$, we have $\eta >0$. Since the loss function is nonnegative,  $l(\cdot)\geq 0$ and $\bar{V}_k\geq 0$ for all $k$.
Taking full expectation and summing from $k=0$ to $k=K-1$ on both sides of
the above inequality, we obtain 
\begin{equation*}
\mathbb{E}[\eta \sum_{k=0}^{K-1}\Vert \tr \ell (\bar{\theta}_{k})\Vert
^{2}]\leq \ell (\bar{\theta}_{0})+L\overline{V}_{0}+\frac{KL{\gamma}^{2}\xi\sigma ^{2}}{N^{2}}(1+\frac{1}{2N}).
\end{equation*}%
We conclude that 
\begin{eqnarray*}
&&\mathbb{E}\left[ \frac{1}{K}\sum_{k=0}^{K-1}\mathbb{E}[\Vert \tr \ell (%
\bar{\theta}_{k})\Vert ^{2}]\right] 
\leq \frac{1}{\eta K}\left[ \ell (\bar{\theta}_{0})+L\overline{V}_{0}+%
\frac{KL\gamma^{2}\xi\sigma ^{2}}{N^{2}}(1+\frac{1}{2N})\right]. \blacksquare
\end{eqnarray*}

\textbf{Proof of Corollary 1}

Since ${\gamma}<\bar{\gamma}_{1}$, it follows that
\al{\k<&2L\mu_{\min}-2a\l_2\mu_{\max}+\frac{2}{3}\Big(2a\l_2\mu_{\max}-L(2\mu_{\min}+\frac{\xi}{2})\Big)\\
	&=\frac{2}{3}L\mu_{\min}-\frac{2}{3}a\l_2\mu_{\max}-\frac{\xi L}{3}<0,
} and from Lemma 2: 
\begin{equation*}
\frac{\left\vert \kappa \right\vert \gamma }{N}\overline{V}_{k}\leq 
\overline{V}_{k}-\mathbb{E}[\overline{V}_{k+1}|\*{\t}_k]+\frac{4\g^{2}\xi}{N}\left\Vert \tr \ell (\bar{\theta}_{k})\right\Vert ^{2}+%
\frac{\gamma^{2}\xi\sigma ^2}{N^{2}}.
\end{equation*}%
Taking full expectation and summing from $k=0$ to $k=K-1$ on both sides of
the above inequality:%
\begin{equation*}
\mathbb{E}\left[ \frac{1}{K}\sum_{k=0}^{K-1}\overline{V}_{k}\right] \leq 
\frac{N}{K|\k|\g}\overline{V}_0+\frac{4\g\xi}{| \kappa | }\Big[\frac{1}{K}\sum_{k=0}^{K-1}\left\Vert \tr \ell (\bar{\theta}_{k})\right\Vert ^{2}+\frac{\sigma ^{2}}{4N}\Big],
\end{equation*}%
and using Theorem 1 we obtain the result. $\blacksquare$

\medskip \textbf{Proof of Proposition 1}


By Assumption 3 and Taylor expansion,
\begin{align*}
\ell (\theta _{k+1})& \leq \ell (\t_k)+\tr \ell (\theta
_{k})\trs(\theta _{k+1}-\t_k)+\frac{L}{2}\left\Vert
\theta _{k+1}-\t_k\right\Vert ^{2} \\
& =\ell (\t_k)-\gamma \tr \ell (\t_k)\trs\sum_{i=1}^{N}\*1\ik g_{i}(\t_k)+\frac{L}{2}\left\Vert \theta
_{k+1}-\t_k\right\Vert ^{2}.
\end{align*}%
Taking conditional expectation on both sides,%
\begin{align*}
\E[\left. \ell (\theta _{k+1})\right\vert {\t}_k]
\leq &\ell(\t_k)-\gamma \tr \ell (\t_k)\trs\sum_{i=1}^{N}\frac{\mu_i}{\mu}\E[\tr \ell(\t_k)+\ve_i]+\frac{L}{2}\E[\left\Vert \theta _{k+1}-\theta_{k}\right\Vert ^{2}| \t_k]\\
=& \ell(\t_k)-\gamma \norm{\tr \ell (\t_k)}^2+\frac{L}{2}\E[\left\Vert \theta _{k+1}-\theta_{k}\right\Vert ^{2}| \t_k].
\end{align*}%
Note that 
\begin{align*}
&\E[ \left\Vert \theta _{k+1}-\t_k\right\Vert^{2}| \t_k]=\E[ \Vert \g\sum_{i=1}^N\*1\ik g_i(\t_k)  \Vert ^{2}| \t_k]\\=&\g^2\sum_{i=1}^{N}\frac{\mu_i}{\mu}\norm{\tr\ell(\t_k)+\ve_i(\t_k)} ^{2} 
\leq\g^2\big(\norm{\tr\ell(\t_k)}^2+\frac{1}{\mu}\sum_{i=1}^{N}\mu_i\s_i^2  \big),
\end{align*}%
it follows that%
\al{\gamma (1- \frac{L\gamma}{2})\Vert \tr \ell(\t_k)\Vert ^{2} 
\leq \ell (\t_k)-\E[\ell(\t_{k+1})|\t_k]+\frac{L\gamma ^{2}}{2\mu}\sum_{i=1}^{N}\mu_i\s_i^2.}
The results follows by taking full expectation and summing from $k=0$ to $%
k=K-1$ on both sides of the above inequality. $\blacksquare $

\vskip 0.2in 

\bibliographystyle{IEEEtran}
\bibliography{AlfredoBIB}

\end{document}